\documentclass[12pt]{amsart}
\usepackage[english]{babel}
\usepackage{amsfonts}
\usepackage{marginnote}
\usepackage{comment}
\usepackage{tabularx}
\usepackage[protrusion=true,expansion=true]{microtype}
\usepackage{enumerate}
\usepackage{bbm}
\usepackage[mathlines,displaymath,pagewise]{lineno}

\usepackage[margin=0.4 cm]{caption}

\usepackage{tikz}
\usepackage{latexsym}
\usepackage[colorlinks=true, pdfstartview=FitV, linkcolor=blue, citecolor=blue, urlcolor=blue,pagebackref=false]{hyperref}

\usepackage{amsmath, amsthm, amssymb, color, graphicx, hyperref, mathrsfs}
\usepackage[margin=1.4in]{geometry}
\DeclareGraphicsExtensions{.jpg,.pdf,.eps}
\graphicspath{{./Figures/}}

\newtheorem{thm}{Theorem}[section]

\newtheorem{prop}[thm]{Proposition}
\newtheorem{lemma}[thm]{Lemma}
 
\newtheorem{corollary}[thm]{Corollary} 

\newtheorem{definition}[thm]{Definition}

\newtheorem{remark}[thm]{Remark}

\newtheorem{claim}[thm]{Claim}

\renewcommand{\ni}{\noindent}

\newcommand{\ave}[1]{{\color{violet}#1}}

\newcommand{\old}[1]{{}}

\newcommand{\C}{\mathbb{C}}

\newcommand{\1}{\mathbf{1}}
\newcommand{\D}{\mathbb{D}}
\newcommand{\E}{\mathbb{E}}
\newcommand{\N}{\mathbb{N}}
\newcommand{\Q}{\mathbb{Q}}

\newcommand{\R}{\mathbb{R}}
\renewcommand{\P}{\mathbb{P}}

\newcommand{\wt}{\widetilde}

\newcommand{\eps}{\varepsilon}
\def\P{\mathbb{P}}
\def\E{\mathbb{E}}
\DeclareMathOperator{\Cov}{Cov}

\def\cE{\mathcal{E}}

\renewcommand{\ni}{\noindent}

\def\@rst #1 #2other{#1}

\newcommand{\aryb}{\begin{eqnarray*}}
\newcommand{\arye}{\end{eqnarray*}}
\def\alb#1\ale{\begin{align*}#1\end{align*}}
\newcommand{\eqb}{\begin{equation}}
\newcommand{\eqe}{\end{equation}}
\newcommand{\eqbn}{\begin{equation*}}
\newcommand{\eqen}{\end{equation*}}

\def\bi{\begin{itemize}}
	\def\ei{\end{itemize}}
\def\bnum{\begin{enumerate}}
	\def\enum{\end{enumerate}}
\def\md{\mid}
\def\Bb#1#2{{\def\md{\bigm| }#1\bigl[#2\bigr]}}

\def\Eb{\Bb\E}
\def \p {{\partial}}

\def\<#1{\langle #1\rangle}

\begin {document}

\title[Negative moments for GMC on fractal sets]{Negative moments for Gaussian multiplicative chaos on fractal sets}

\date{}

\author[Garban]{Christophe Garban}
\address[Christophe Garban]{Univ Lyon, Universit\'e Claude Bernard Lyon 1, CNRS UMR 5208, Institut Camille Jordan, 69622 Villeurbanne, France}
%
\author[Holden]{Nina Holden}
\address[Nina Holden]{Massachusetts Institute of Technology, 77 Massachusetts Avenue, Cambridge, MA 02139-4307, USA}
%
\author[Sep\'ulveda]{Avelio Sep\'ulveda}
\address[Avelio Sep\'ulveda]{Univ Lyon, Universit\'e Claude Bernard Lyon 1, CNRS UMR 5208, Institut Camille Jordan, 69622 Villeurbanne, France}
%
\author[Sun]{Xin Sun}
\address[Xin Sun]{Columbia University, 2990 Broadway, New York, NY, 10027, USA.}

\maketitle

\begin{abstract}
    The objective of this note is to study the probability that the total mass of a  sub-critical Gaussian multiplicative chaos (GMC) with arbitrary base measure $\sigma$ is small.  When $\sigma$ has some continuous density w.r.t Lebesgue measure, a scaling argument shows that the logarithm of the total GMC mass is sub-Gaussian near $-\infty$.  However, when $\sigma$ has no scaling properties, the situation is much less clear. In this paper, we prove that for any base measure $\sigma$, the total GMC mass has negative moments of all orders. 

\end{abstract}

\section{Introduction}
\subsection{Context}
Gaussian multiplicative chaos (GMC) measures are a way to give meaning to the exponential of random ``generalised functions'' that cannot be defined pointwise. They were introduced by Kahane in \cite{KAH}, and in recent times they have seen a revived interest with two dimensional Liouville quantum gravity \cite{DS}.
	
	GMC measures can be informally expressed as ``$\mu_{\gamma h}^\sigma:=e^{\gamma h(x)}\sigma(dx)$'', where $h$ is a log-correlated Gaussian field and $\sigma$ is a finite measure. One of the first questions addressed in the study of GMC theory concerns the triviality of $\mu_{\gamma h}^\sigma$. In \cite{KAH}, it was shown that if there exists $d$ such that $\mathcal{E}_d(\sigma):=\iint\|x-y\|^d\sigma(dx)\sigma(dy)<\infty$, then almost surely $\mu_{\gamma h}^\sigma\neq0$ if $\gamma<\sqrt{2d}$.

There have already been notable instances where the non-triviality of the GMC measure has been quantified. To do this, one needs to study the tail behaviour of the total mass of the GMC near 0. In most of these cases, the base measure $\sigma$ has been taken as the Lebesgue measure, $\lambda$, restricted to some open set $D$. Below is a non-exhaustive list of important results on the tails of GMC.
\begin{itemize}
	\item First, the existence of all negative moments for Gaussian multiplicative chaos $\mu^\lambda_{\gamma h}(\D)$ has been proved in \cite{RobertVargas}.

	\item In \cite{DS}, an important step in the proof of the {\em KPZ relation} is the following result (Lemma 4.5 in \cite{DS}): for any open ball $B\subset \D$, there exists $c=c_{\gamma,B}>0$ such that for any $\eps\in(0,1)$, $\P(\mu_{\gamma h}^\lambda(B)<\eps) \leq e^{- c (\log \frac 1 \eps)^2}$. Here, $h$ is a zero-boundary GFF in $\D$. See also Lemma 4.3 of \cite{Anotes} for a review.

	\item In \cite{Remy}, the exact density of the total mass of the GMC (for a GFF with Neumann boundary conditions) on the unit circle $\p \D$ is obtained, thus answering a conjecture by Fyodorov and Bouchaud \cite{fb08}. In particular, it shows that if $Y_\gamma$ is the total mass of the GMC on $\p \D$ with parameter $\gamma/ 2$, then $\P(Y_\gamma  < \eps) \approx \exp(-c_\gamma \eps^{-4/\gamma^2})$.

\item In the opposite direction, the upper tails (i.e.\ $\P(\mu_{\gamma h}^\lambda(\D) \geq x)$ for large $x$) have been studied in detail in the recent \cite{RVtails}. 
\end{itemize}

\subsection{Results}
\begin{figure}
	\begin{center}
		\includegraphics[width=0.75\textwidth]{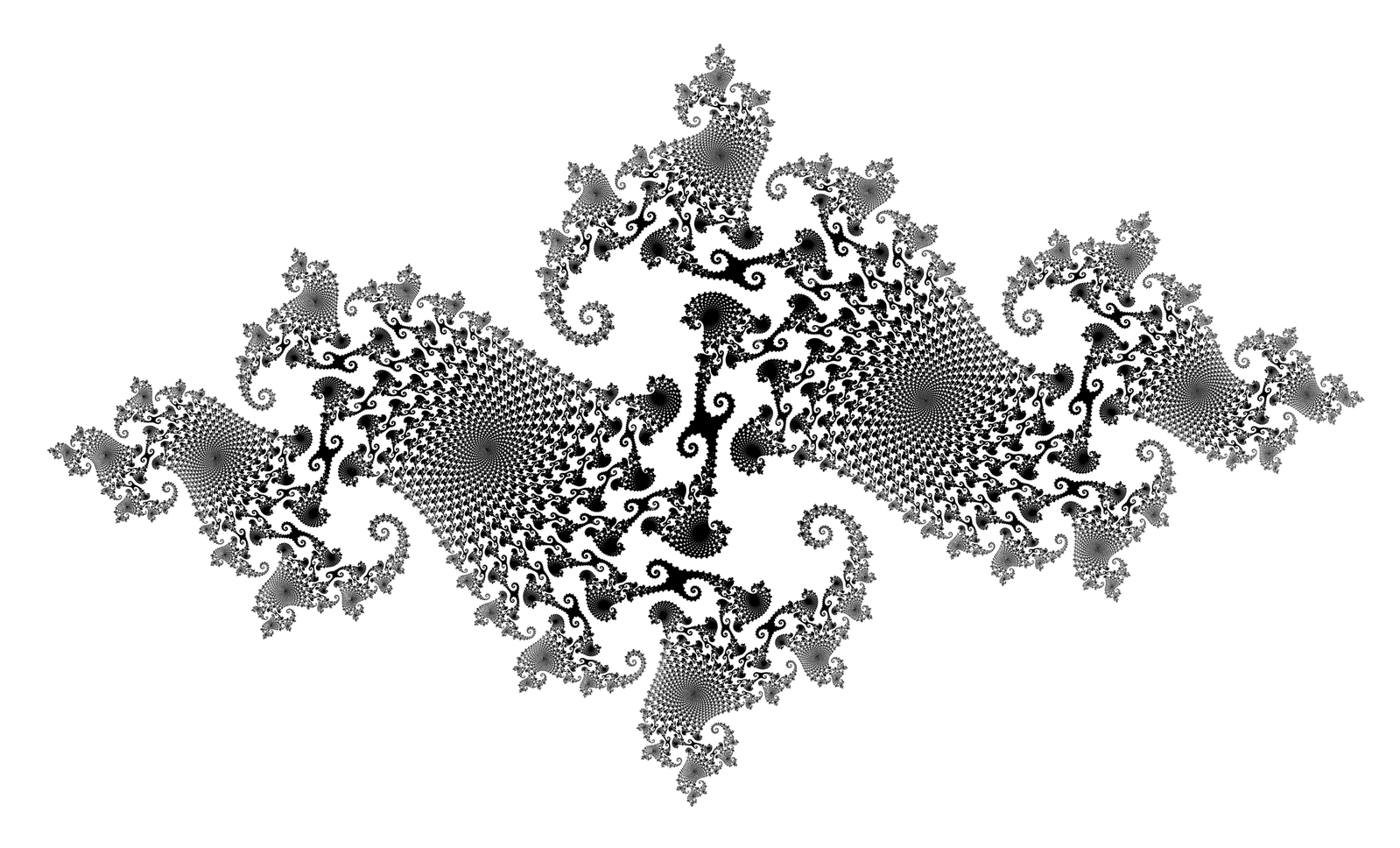}
	\end{center}
	\caption{To illustrate our main result: suppose you are given a Borel measure $\sigma$ supported on the above Julia set $J$ whose Hausdorff dimension is close to two. Suppose this measure is well-spread on $J$ in the sense that $\cE_d(\sigma)<\infty$ for some dimension $d<2$. If $h$ is a 2d GFF, this allows us to consider a Gaussian multiplicative chaos measure $\mu_{\gamma h}=\mu_{\gamma h}^\sigma$ supported on $J$ if $\gamma<\sqrt{2d}$. As  this Julia set does not have any {\em exact} scaling invariance, standard scaling tools fail to study the tails of $\P(\int_{\C} e^{\gamma h(x)} \sigma(dx) \leq \eps)$ when $\eps \to 0$. In this work, we obtain quantitative (sub-optimal) bounds on these tails for any base measure $\sigma$  satisfying only $\sigma(\D) <\infty$ and $\cE_d(\sigma)<\infty$. The application we have in mind is not to consider GMC measures on Julia sets but rather GMC measures on {\em spectral samples} of critical percolation which were introduced in \cite{GPS} (see \cite{LDP}). As these spectral samples are typically fractal sets whose scaling properties are poorly understood, the investigation in this work cannot be avoided. \copyright  Prokofiev \href{https://creativecommons.org/licenses/by-sa/3.0}{CC BY-SA 3.0}.}
	\label{f.julia}
\end{figure}

All the presented results concerning the tails of GMC measures rely on techniques which require some forms of exact scaling invariance of GMC measures. The goal of this paper is to obtain bounds on GMC measures which are supported on fractal sets without any a priori invariance under scaling. This problem arose in \cite{LDP}, where we need tail estimates for GMC measures which are defined on the {\em spectral samples} of critical percolation \cite{GPS}. These spectral samples are random fractal sets of dimension $3/4$ which do not have any {\em quenched} scaling properties. An entirely different analysis, which is not built out of re-scaling arguments, is thus required  (see also Figure \ref{f.julia} for an illustration of our main result/motivations).

For simplicity, we work in most of this paper with the (zero-boundary) Gaussian free field (GFF) $h$ in the unit disk $\D$, i.e., the centred Gaussian random distribution $h$ in $H^{-\delta}(\D)$ which has covariance kernel
\begin{align*}\label{}
\Cov(h(x),h(y)) = G_\D(x,y) = \log \left|\frac {1-x\bar y} {x-y}\right |.
\end{align*}
Here $\delta>0$ and $H^{-\delta}(\D)$ is the Sobolev space with index $-\delta$. See \cite{She} for an introduction to the GFF. If needed, the results can be easily extended to any log-correlated field (see Remark \ref{r.ext}). 
Our setup is as follows: we consider a fixed Borel measure $\sigma$ with compact support in $\D$ and such that it is $d$-dimensional in the sense that $\cE_d(\sigma)<\infty$ (see Definition \ref{d.energy}). It follows from the works of \cite{KAH,DS} that if $\gamma<\sqrt{2d}$ one can define the GMC measure $\mu_{\gamma h}^\sigma := ``e^{\gamma h(x)} \sigma(dx)"$. One way to make proper sense of the exponential of this distribution is through a limiting procedure. Let $h$ be a zero-boundary GFF in $\D$ and define $h_\epsilon(z)=h*\rho_\epsilon$, where $\rho_\epsilon$ is a smooth mollifier. Then, for any test function $f\in C_c(\D)$,
\begin{align*}\label{}
\int_\D f(x) e^{\gamma h_\eps(x) - \frac {\gamma^2} 2 \Eb{h_\eps^2(x)}} \sigma(dx) \xrightarrow{\eps \to 0} \int_D f(x) \mu_{\gamma h}^\sigma(dx),
\end{align*}
where the limit holds a.s.\ and in $L^1$ as $\eps \to 0$. See the references \cite{DS,RVreview,Nath,Anotes} for general background on Gaussian multiplicative chaos (GMC) as well as other possible regularisation procedures. 

\medskip 

As stated before, the main result of this paper gives an estimate on the tails of the GMC measure near 0. 
\begin{thm}\label{thm:: laplace}
	Let $\sigma$ be a $d$-dimensional measure (i.e.\ such that $\cE_d(\sigma)<\infty$) with compact support in $ \D$, let $\gamma<\sqrt{2d}$, let $h$  be a (zero-boundary) GFF in $\D$, and let $\mu_{\gamma h}^\sigma$ be the GMC measure with parameter $\gamma$ associated to $h$ and base measure $\sigma$. Then, for all $n\in \N$ there exists $t_0=t_0(\sigma,\gamma,n)>0$ such that for all $t\geq t_0$
	\begin{equation*}
	\E\left[e^{-t\mu_{\gamma h}^\sigma(\D)} \right] \leq t^{-n}.
	\end{equation*}
	In particular, $1/\mu_{\gamma h}^\sigma(\D)$ has moments of all orders.
\end{thm}

In some cases, the dependence of $t_0$ as a function of the parameters $\sigma,\gamma,n$ can be made quantitative: see for example Lemma \ref{lemma:: quantitative} or Corollary \ref{e.simplified}.
\begin{remark}\label{r.1}
Using Markov's inequality, the theorem implies readily that for any $n\geq 1$, there exists $\eps_0\in (0,1)$ such that for any $\eps<\eps_0$, $\P(\mu_{\gamma h}^\sigma (\D)< \eps) \leq \eps^n$. 
\end{remark}

Finally, let us note that our setup can be generalised to any log-correlated field.
\begin{remark}\label{r.ext}
	Using Kahane's convexity inequality (Proposition \ref{Prop:: Kahane}), it is straightforward to extend it to Neumann GFF in $\D$ or to any GFF in other simply connected domains $D \subset\C$. Note also that for zero-boundary GFF, the assumption that $\sigma$ has compact support in $\D$ can be removed by an easy dichotomy argument (separating the possible mass on $\p \D$ from the mass in its interior).  
\end{remark}

\ni
\textbf{Acknowledgments:}
The research of C.G.\ is supported by the 
ANR grant \textsc{Liouville} ANR-15-CE40-0013 and the ERC grant LiKo 676999.
The research of N.H.\ is partially supported by a fellowship from the Norwegian Research Council.
The research of A.S.\ is supported by the ERC grant LiKo 676999.
The research of X.S.\ is supported by Simons Society of Fellows under Award 527901.
The work on this paper started during the visit of N.H.\ and X.S.\ to Lyon in November 2017. They thank for the hospitality and for the funding through the ERC grant LiKo 676999. 

\section{Preliminaries}

\subsection{Energy and $d$-dimensional measures}

\begin{definition}\label{d.energy}
For any Borel measure $\sigma$ in $\D$ and any $d>0$, we define its \textbf{$d$-energy} to be 
\begin{align*}
\cE_d(\sigma):=\iint \frac 1 {\| x- y\|_2^{d}} \sigma(dx) \sigma(dy).
\end{align*}
We will say that a measure $\sigma$ on $\D$ is \textbf{$d$-dimensional} if $\cE_d(\sigma)<\infty$. 
\end{definition}

We will rely extensively in this work on the following notion of local energy.
\begin{definition}[Local $d$-energy]\label{d.local}
For any base point $\bar x \in \D$, any Borel measure $\mu$ on $\D$ and any real
$\beta>0$, we define  
\begin{align*}\label{}
\phi_\beta(\bar x,\mu):= \int_\D \frac 1 {\| \bar x-  x\|_2^{\beta}}  \mu(dx).
\end{align*}
\end{definition}

\subsection{Gaussian free field} 
We say that $h$ is a GFF in a domain $D\subseteq \R^2$, if it is a centred Gaussian ``generalised function'' such that for any smooth function $f$ \[\E[(h,f)^2]=\iint f(x)G_D(x,y)f(y) dxdy.\] Here $G_D(x,y)$ is the Green's function with zero-boundary in $D$ with normalisation such that $G_D(x,y) \sim |\log(|x-y|)|$ as $y\to x$. Furthermore, if either $x$ or $y$ do not belong to $D$ then $G_D(x,y)=0$. Let us note that if $D\subseteq D'$ then $G_D(x,y)\leq G_{D'}(x,y)$. 

The main property of the GFF we use in this paper is its Markov property.
\begin{lemma}[Markov property]\label{lemma:: Markov} Let $h$ be a GFF in $D$ and let $A\subseteq \overline D$ be a closed (deterministic) set. Then, the restriction of $h$ to $D\backslash A$ can be written as the independent sum of $h^A$ and $h_A$ where $h^A$ has the law of a GFF in $D\backslash A$ and $h_A$ is a harmonic function (thus continuous) in $D\backslash A$.
\end{lemma}


\subsection{Three useful inequalities}

The first inequality we shall rely on is the following famous inequality for general centred Gaussian processes:
\begin{thm}[Borell-TIS inequality, see \cite{Adler}]\label{th.TIS}
Let $(X_{n})_{n\in \N}$ be a centred Gaussian field. If a.s.\ $\sup X_{n} <\infty$, then $\E\left[\sup X_n\right] <\infty$ and for any $t>0$
	\begin{equation*}
	\E\left[e^{t \sup X_n} \right]\leq e^{t \E\left[ \sup X_n\right] + \frac{t^2\tilde \sigma^2}{2}},
	\end{equation*}
	where $\tilde \sigma^2 := \sup_n \text {Var}\left[X_n^2\right]<\infty$.
\end{thm}

\ni
The second inequality is the so-called FKG inequality proved by Pitt in 1982 \cite{Pitt}.
\begin{thm}[FKG inequality, \cite{Pitt}]\label{thm:: FKG}
 Let $(X_\lambda)_{\lambda \in \Lambda}$ be a centred Gaussian field with $\E\left[X_\lambda X_{\lambda'}\right]\geq 0$ for all $\lambda,\lambda' \in \Lambda$. Then, if  $f,g$ are two bounded, increasing measurable functions,
	\begin{equation*}
	\E\left[f((X_{\lambda})_{\lambda \in \Lambda})g((X_{\lambda})_{\lambda \in \Lambda})\right]\geq \E\left[f((X_{\lambda})_{\lambda \in \Lambda})\right] \E\left[g((X_{\lambda})_{\lambda \in \Lambda})\right].
	\end{equation*}
\end{thm}

\ni
Finally, our last inequality will be the following result proved by Kahane in 1982 \cite{KAH}. See also \cite{RVreview} or Proposition 6.1 of \cite{Anotes}.
	\begin{prop}[Kahane's convexity inequality]\label{Prop:: Kahane} Let $(h^1_{\lambda})_{\lambda \in \Lambda}$ and $(h^2_\lambda)_{\lambda\in \Lambda}$ be two log-correlated centred Gaussian  fields such that their covariance kernels satisfy $\mathcal C_1(x,y)\leq \mathcal C_2(x,y)$. Then, for any convex function $f$, we have that
		\begin{equation*}
		\E\left[f\left (\mu_{h^{1}_{\lambda}}(\D)\right ) \right] \leq \E\left[f\left (\mu_{h^{2}_{\lambda}}(\D)\right ) \right]. 
		\end{equation*}
	\end{prop}

\subsection{A useful change of measure} A key ingredient in our proof will be a change of measure associated to the total GMC mass.

Let $\P$ be a measure where $h$ is a GFF in $D$. Let us describe how the law of $h$ changes when one weights $\P$ by the total GMC mass $\mu^\sigma_{\gamma' h}(\D)$ of parameter $\gamma'<\sqrt{2d}$ (not necessarily equal to $\gamma$ for our later applications of the change of measure below). The following result is Theorem 17 of \cite{Shamov} or Proposition 3.1 of \cite{Anotes}. It goes back to the work by Kahane and Peyri\`ere \cite{KP}.

\begin{prop}\label{law Q} 
Let $\sigma$ be any $d$-dimensional measure. For any  $\gamma' <\sqrt{2d}$, we introduce the following probability measure on $H^{-1}(\D)$, 
	\begin{linenomath*}
		\begin{align*}
		\frac{d\Q_{\gamma'}}{d\P}(h) := \frac {\mu^\sigma_{\gamma' h}(\D)} {\sigma(\D)},
		\end{align*}
	\end{linenomath*}
where $\P$ is the law of an unbiased GFF in $\D$ with Dirichlet boundary conditions. (N.B. For any $\gamma' < \sqrt{2d}$, the probability measure $\Q_{\gamma'}$ is well defined and absolutely continous w.r.t $\P$).  
Then, under the new measure $\Q_{\gamma'}$, if $h \sim \Q_{\gamma'}$, we have the identity in law,
	\begin{linenomath*}
		\begin{align*}
h(\cdot)\overset{(law)}=\hat h(\cdot) + \gamma'\, G_{\D}(\bar x,\cdot),
		\end{align*}	\end{linenomath*}
		where, under the law $\Q_{\gamma'}$, $\hat h$ is an (unbiased) GFF in $\D$ and $\bar x$ is a random point independent of $\hat h$ and sampled according to $\bar x\sim \sigma(dx)/\sigma(\D)$.
\end{prop}

Let us state an important consequence of this change of measure.
\begin{lemma}\label{obvious}
Let $\sigma$ be a $d$-dimensional measure ($\cE_d(\sigma)<\infty$) in $\D$. 
For any $\gamma,\gamma' < \sqrt{2d}$, 
if $\hat h$ is a (zero boundary) GFF in $\D$, $\bar x \sim \sigma(dx)/\sigma(\D)$ is a random point independent of $\hat h$, and $h(\cdot):= \hat h(\cdot) + \gamma'\, G_{\D}(\bar x, \cdot)$, then the following identity holds a.s.
\begin{align}\label{e.obvious}
\int_{\D} e^{\gamma \gamma' G_{\D}(\bar x, x)} \mu_{\gamma \hat h}^\sigma(dx) = \mu_{\gamma h}^\sigma(\D).
\end{align}
Additionally, if we define $\bar \beta:=\max\{\sqrt{2d}\gamma, d\}$ we have for any $\beta<\bar \beta$, 
\begin{align}\label{e.beta}
\int_{\D} e^{\beta G_{\D}(\bar x, x)} \mu_{\gamma \hat h}^\sigma(dx) < \infty \text{  a.s.}
\end{align}
\end{lemma}
\begin{remark}\label{}
Note perhaps surprisingly that for some values of $\gamma$, the presence of the multiplicative chaos allows us to integrate more singular kernels $1/\|\bar x - x\|^\beta$ than what is expected readily from the energy bound. 
\end{remark}
On the formal level, the identity~\eqref{e.obvious} looks obvious, but its (short) proof is normally omitted. We thus include a short justification below.
\medskip

\ni
{\em Proof.}
From Proposition \ref{law Q}, we have that the random field $h$ is absolutely continuous w.r.t an (unbiased) zero-boundary GFF on $\D$. In particular, we have that the measure $\mu_{\gamma h, \eps}^\sigma(dx):= e^{\gamma h_\eps(x) - \frac {\gamma^2} 2 \Eb{\hat h_\eps(x)^2}} \sigma(dx)$ converges in probability to $\mu_{\gamma h}^\sigma$, where $h_\epsilon=h*\rho_\epsilon$ for $\rho_\epsilon$ a smooth mollifier. 
For any $\delta>0$, this implies that a.s. 
\begin{align*}\label{}
\mu_{\gamma h, \eps}^\sigma(\D \setminus B(\bar x, \delta)) \xrightarrow{\eps \to 0} \mu_{\gamma h}^\sigma(\D \setminus B(\bar x, \delta)).
\end{align*}
Now, we may rewrite $\mu^\sigma_{\gamma h, \eps}(\D \setminus B(\bar x, \delta))$ as follows
\begin{align*}\label{}
\mu^\sigma_{\gamma h, \eps}(\D \setminus B(\bar x, \delta)) 
& = \int_{B(\bar x, \delta)^c} e^{\gamma \gamma' G_\D^{\eps}(\bar x, x)} \mu^\sigma_{\gamma \hat h, \eps}(dx),
\end{align*}
where $G_\D^\eps(\bar x, \cdot):= \rho_\eps * G_\D(\bar x, \cdot)$. 
By the convergence in probability of $\mu^\sigma_{\gamma \hat h,\eps} \to \mu^\sigma_{\gamma \hat h}$ and the absence of singularity of the exponential term in $B(\bar x, \delta)^c$, we obtain by taking $\eps  \to 0$ the identity
\begin{align*}\label{}
\mu^\sigma_{\gamma h}(\D \setminus B(\bar x, \delta))  = \int_{B(\bar x, \delta)^c} e^{\gamma \gamma' G_\D(\bar x, x)} \mu^\sigma_{\gamma \hat h}(dx).
\end{align*}
Now, we conclude the proof of the identity~\eqref{e.obvious} by letting $\delta\to 0$ using the monotone convergence theorem together with the a.s.\ absence of Dirac point masses for both $\mu^\sigma_{\gamma h}$ and $\mu^\sigma_{\gamma \hat h}$ ($\gamma<\sqrt{2d}$).

	When $\beta <\sqrt{2d}\gamma$, the second identity~\eqref{e.beta} follows from~\eqref{e.obvious} plus the fact that $\mu_{\gamma h}(\D) <\infty$ a.s. When $\beta < d$ (in fact $\beta$ may even  be equal to $d$ here), we note that
\begin{align*}\label{}
\E\left[ \int_{\D} e^{\beta G_{\D}(\bar x, x)} \mu_{\gamma \hat h}(dx)\right]
&= \frac 1 {\sigma(\D)} \iint e^{\beta G_\D(x,y)} \sigma(dx) \sigma(dy) \\
&\leq \frac {2^2} {\sigma(\D)} \iint \frac 1 {\|x -y\|^\beta} \sigma(dx) \sigma(dy) <\infty,
\end{align*}
where we used the fact that $G_\D(x,y) \leq \log(2) + |(\log |x-y|)|$ together with the fact that $\beta \leq d\leq 2$. \qed

\section{GMC measures have negative moments of some order}
The goal of this section is to prove that there exists an (explicit) exponent $\eta>0$ such that the Laplace transform of $\mu_{\gamma h}^\sigma(\D)$ is $O(t^{-\eta})$. Even though it is not necessary for the proof of Theorem \ref{thm:: laplace}, we are going to be quantitative. This will be important in particular to obtain quantitative ergodic bounds in our coming work  {\em Liouville dynamical percolation} \cite{LDP}. Moreover, it is a key new input in the upcoming revised version of \cite{BSS} proving that the GMC measure is determined by the GFF under mild conditions.
As such, the lemma below is of independent interest. Furthermore, let us remark that the proof of this lemma uses the classical change of measure of Proposition \ref{law Q} in a way which  to our knowledge is new.

 In order to state the result, let us recall that for any $0<\gamma<\sqrt{2d}$, we defined in Lemma \ref{obvious},  $\bar \beta:=\max\{d, \sqrt{2d} \gamma \}$.

\begin{lemma}\label{lemma:: quantitative}
Let $\sigma$ be a Borel measure with compact support in $\D$ such that  $\sigma(\D)<\infty$ and $\cE_d(\sigma)<\infty$ for some $d\leq 2$. 
First, fix any choice of $\delta>0$ and $\beta\in (\gamma^2,\bar \beta)$ and define the following exponent 
\[
\eta=\eta_{\delta,\beta}:=\frac {\beta-\gamma^2} {\beta + \gamma^2 \delta}.
\]
Let $h$ be a zero-boundary GFF in $D\subseteq \D$ (recall that $h$ is 0 outside of $D$) and let $\mu=\mu_{\gamma h}^\sigma$ be the GMC measure of parameter $\gamma$ and base measure $\sigma$. Then, there exists $t_0>0$ such that for any 
$t \geq t_0$, 
\begin{linenomath*}
	\begin{align*}
	\E\left[\exp\left (-t\mu(\D)\right ) \right] \leq \frac{2^5}{\sigma(\D)t^{\eta}}.
	\end{align*}
\end{linenomath*}
Furthermore, one can take $ t_0= 2^4 s_0^{1/\eta}$ where $s_0$ is a positive real number so that  
\begin{linenomath*}
\begin{align}\label{e.nasty}
\P\left (\phi_{\beta}(\bar x,\mu)\leq \frac 1 {2^4} s_0^{\delta} \right )\geq 1/2.
\end{align}
\end{linenomath*}
(N.B.  Recall the definition of $\phi_\beta$ in Definition \ref{d.local}. The existence of a positive $s_0$ satisfying the above condition is ensured by Lemma \ref{obvious}.) 	
\end{lemma}

The advantage of this lemma as compared to our main result (Theorem \ref{thm:: laplace}) is that it quantifies the condition on $t_0=t_0(\sigma,\gamma, \eta)$. However, the exponent obtained is not very good (it is $<1$) and the condition~\eqref{e.nasty} behind the definition of $t_0$ is hard to digest. Let us then state the following corollary of the proof of Lemma \ref{lemma:: quantitative} in the $L^2$-regime $\gamma<\sqrt{d}$, where the $t_0$-dependence becomes much more readable.

\begin{corollary}[Simplified quantitative estimate in the $L^2$ regime $\gamma<\sqrt{d}$]\label{e.simplified}

Let $\sigma$ be a Borel measure with compact support in $\D$ such that $\sigma(\D)<\infty$ and $\cE_d(\sigma)<\infty$ for some $d\leq 2$. Then, for any $\gamma<\sqrt{d}$, if 
\[
\eta:= \frac{d-\gamma^2}{d+\gamma^2},
\]
we have 
\begin{align*}
	\E\left[e^{-t\mu_{\gamma h}^\sigma(\D)} \right] \leq \frac{2^5}{\sigma(\D)t^{\eta}},
	\end{align*}
for any
\[
t\geq t_0:= 2^4 \left [ 2^{5} \frac{\cE_d(\sigma)}{\sigma(\D)}\right ]^{1/\eta}.
\]
\end{corollary}

\ni
{\em Proof of Lemma \ref{lemma:: quantitative}.} As $G_D(x,y)\leq G_\D(x,y)$ we can use Kahane's convexity inequality (Proposition \ref{Prop:: Kahane}) to reduce ourselves to the case where $h$ is a GFF in $\D$. Define $\Q:=\Q_\gamma$ as in Proposition \ref{law Q}. Using the fact that $xe^{-x s}\leq e^{-1}/s$ for any $x\geq 0, s>0$, we have (with $\mu=\mu_{\gamma h}^\sigma$), 
	\begin{linenomath*}
		\begin{align*}
		\Q\left[\exp(- s \mu(\D))\right]&=\frac{ \E\left[\mu(\D)\exp(- s \mu(\D)) \right]}{\sigma(\D)}\leq \frac{e^{-1}}{\sigma(\D) s}.	
		\end{align*}
	\end{linenomath*}
	Thus, thanks to the identity \eqref{e.obvious} we obtain for any $s>0$ the bound
	\begin{linenomath*}
		\begin{align}\label{change measure result laplace}
\E\left [\exp\left (-s \int_\D e^{\gamma^2 G_{\D}(\bar x, x)} \mu(dx)\right ) \right] \leq \frac{e^{-1}}{\sigma(\D)s}.
		\end{align}	\end{linenomath*}
	
This is almost what we wish to prove except the log-singularity at $\bar x$ plays against us. Indeed, it may have the effect that $\E\left [\exp\left (-s\int_\D e^{\gamma^2 G_{\D}(\bar x, x)} \mu(dx)\right ) \right]$ is much smaller than  $\E\left [\exp\left (-s \mu(\D) \right ) \right]$. To analyse the impact of this log-singularity at $\bar x$, take $r>0$ to be chosen later and let us  introduce the following event:
\begin{align*}\label{}
C(s,r) : = \left\{\int_{B(\bar x,r)} \exp(\gamma^2 G_{\D}(\bar x,x)) \mu(dx)\leq s^{\delta} r^{\beta-\gamma^2} \right\},
\end{align*}
i.e., the event that  $\mu=\mu_{\gamma h}^\sigma$ does not put a lot of mass in $B(\bar x, r)$. (Here $h \sim \P$ is a non-biased GFF with zero boundary condition). Thanks to the fact that $G_{\D}(\bar x, x)\leq \log(\|\bar x-x\|)+\log(2)$,  we get the following upper bound on the event $C(s,r)$, 
	\begin{linenomath*}
		\begin{align*}
		\int_{\D} e^{\gamma^2 G_{\D}(\bar x, x)} \mu(dx)&\leq 2^{\gamma^2}\int_{\D\backslash B(\bar x, r)} |\bar x-x|^{-\gamma^2}\mu(dx)+s^\delta r^{ \beta-\gamma^2}\\
		&\leq 2^4r^{-\gamma^2}\mu(\D)+s^\delta r^{\beta-\gamma^2}.		\end{align*}
	\end{linenomath*}
	Now, it makes sense at this stage to tune $r$ in a way such that $s\cdot s^{\delta}r^{ \beta-\gamma^2}=1$. In other words, let $r:=s^{-L}$, with $L:=(1+\delta)/(\beta-\gamma^2)$. By doing this and inserting it into \eqref{change measure result laplace} we obtain for any $s>0$, 
	\begin{linenomath*}
		\begin{align*}
		\E\left[\exp(-2^4 s^{1+\gamma^2 L}\mu(\D))\1_{C(s,s^{-L})} \right] \leq \frac{1}{\sigma(\D)s}.
		\end{align*}
	\end{linenomath*}
As $h$ is the zero-boundary GFF in $\D$, it has pointwise positive correlations and thus satisfies the FKG inequality  (Theorem \ref{thm:: FKG}). Since both $\1_{C(s,s^{-L})}$ and $-\mu_{\gamma h}^\sigma(\D)$ are decreasing functions of the field $h$, we have 
	\begin{linenomath*}
		\begin{align}\label{e.FKG}
		\E\left[\exp(-2^4 s^{1+\gamma^2L}\mu(\D))\right] & \leq \frac{\E\left[\exp(-2^4 s^{1+\gamma^2L}\mu(\D))\1_{C(s,s^{-L})} \right]}{\P(C(s,s^{-L}))} \nonumber\\
		& \leq \frac{1}{\sigma(\D)\P(C(s,s^{-L}))s}.
		\end{align}
	\end{linenomath*}
	
We face a potential difficulty here: the function $s\mapsto \P(C(s,s^{-L}))$ does not have any obvious monotonicity. Yet, we shall argue below that it is bounded from below by the following monotone function of  $s>0$: 	
	\begin{linenomath*}
		\begin{align*}
		\P(C(s,s^{-L}))\geq \P\left (\phi_{\beta}(\bar x,\mu) \leq \frac 1 {2^4} s^{\delta} \right ).
		\end{align*}
			\end{linenomath*}
To prove this inequality, suppose the event $\phi_{\beta}(\bar x,\mu_{\gamma h})\leq 2^{-4}s^{\delta}$ occurs. This implies that for any radius $r\in(0,1)$, 
\begin{align*}\label{}
\int_{B(\bar x, r)} e^{\gamma^2 G_\D(\bar x, x)} \mu(dx) 
& \leq 2^{\gamma^2} \int_{B(\bar x, r)} \frac 1 {\|\bar x - x\|^{\gamma^2}} \mu(dx)  \\
& \leq 2^4 r^{\beta-\gamma^2} \int_\D  \frac 1 {\|\bar x - x\|^{\beta}} \mu(dx) \\
& = 2^4 r^{\beta -\gamma^2} \phi_\beta(\bar x, \mu)  \leq  r^{\beta - \gamma^2} s^\delta.
\end{align*} 		
In particular for all $r>0$, $C(s,r)$ is satisfied once $\phi_{\beta}(\bar x,\mu) \leq 2^{-4}s^{\delta}$ holds. We see that for all $s\geq s_0$, we have from the above domination and the definition of $s_0$ \eqref{e.nasty}, that $\P(C(s,s^{-L})) \geq 1/2$ . By noticing that $1+\gamma^2 L = 1/\eta$ and using the change of variables, 
\begin{align*}\label{}
 t=2^4 s^{1/\eta}   \;\;\;\; \text{     and      } \;\;\;\; t_0=2^4 s_0^{1/\eta}  
\end{align*}
 in \eqref{e.FKG}. We obtain that for any $t\geq t_0$, 
\begin{align*}\label{}
\E\left[\exp(-t \mu(\D))\right] & \leq \frac{2}{\sigma(\D) (t/2^4)^\eta} \leq \frac{2^5} {\sigma(\D) t^\eta},
\end{align*}
which concludes our proof.  \qed


\medskip
\ni
{\em Proof of Corollary \ref{e.simplified}.}
We will rely on the above proof and set the parameters as follows. Let $\delta:=1$ and $\beta:=d$ (note that when $\gamma\leq \sqrt{d/2}$, $\beta=\bar \beta = d$ was in fact not allowed in the previous proof, but in the present $L^2$-regime $\gamma<\sqrt{d}$, it will turn out to be harmless). Following the exact same proof, it only remains to check that if $\gamma < \sqrt{d}$, then for any $s\geq  s_0:= 2^{5} \frac{\cE_d(\sigma)}{\sigma(\D)}$, one has 
\[
\P\left (\phi_{d}(\bar x,\mu) \leq \frac s {2^4}  \right )\geq 1/2.
\]
Indeed, it follows directly from Markov's inequality that $\P\left (\phi_{d}(\bar x,\mu)  > 2^{-4}s_0  \right )$ is smaller than or equal to
\begin{align*}\label{}
\frac{2^4} {s_0}
\Eb{\phi_d(\bar x, \mu)} =  \frac{2^4} {s_0\sigma(\D)} \int_\D \sigma(d\bar x) \E\left[ \int_\D \frac 1 {\| \bar x-  x\|_2^{d}}  \mu(dx)\right]  =  \frac{2^4\cE_d(\sigma)} {s_0\sigma(\D)} = 1/2.
\end{align*}
This concludes the proof of the corollary. \qed

\section{Proof of the main result}
In this section, we prove Theorem \ref{thm:: laplace}. We will use a bootstrap argument building on Lemma \ref{lemma:: quantitative}. This part of the proof is close to the classical setting where $\sigma$ has some continuous density w.r.t Lebesgue measure. 

\medskip
\ni
{\em Proof of Theorem \ref{thm:: laplace}.}
Fix $\delta>0$ and $\beta \in (\gamma^2, \bar \beta)$. Let  $\eta:= ( \beta-\gamma^2)/(\beta +\gamma^2\delta)$ as   in Lemma \ref{lemma:: quantitative}. Let us show by induction that for all $d$-dimensional measures $\sigma$ and all $n\in \N$, there exists $K=K(\sigma,n)>0$ such that for all $t\geq 0$ and all $h$ GFF in $D\subseteq \D$,
\begin{equation}\label{eq:: induction}
	\E\left[e^{-t\mu_{\gamma h}^\sigma(\D)} \right] \leq \frac{K}{t^{ 2^n\eta}},
\end{equation}
	When $n=0$, \eqref{eq:: induction} follows from Lemma \ref{lemma:: quantitative} and Kahane's convexity  inequality (Proposition \ref{Prop:: Kahane}).
	
	Let us assume \eqref{eq:: induction} is true for $n\in \N$, WLoG we can assume that $h$ is a GFF in $\D$. Now, for any $\zeta>0$ and any half-plane $H\subset \C$, define $H^+_\zeta$ (resp. $H^-_{\zeta}$) as the set of points in $H$ (resp. in $\C \setminus H$) that are at distance at least $\zeta$ from $\partial H$. 
	\old{We also write $H$ and $D\setminus H$ as $H^+$ and $H^-$ respectively.}
	For simplicity, let us first assume the following claim:	
	\begin{claim}\label{claim::separation}If $\sigma$ is a Borel measure in $\D$ with no atoms, there exists $\zeta>0$ and some half-plane $H\subseteq \C$ such that $4\sigma(H^\pm_\zeta \cap \D)>\sigma(\D)$.
	\end{claim}
	
	Take $\zeta>0$ and $H$ as in the claim. Thanks to Lemma \ref{lemma:: Markov}, we can write $h=h^{+} + h^{-}+h_{\partial H}$, where all terms are independent, 
$h^{+}$ (resp. $h^-$) is a zero-boundary GFF in $H\cap \D$ (resp. $\D\setminus H$), and $h_{\partial H}$ is a harmonic function in $\D \backslash \partial H$. 
	Assuming the mollifier for $h$ is the circle average and using the fact that $h_{\p H}$ is harmonic, we have 
	\[\E\left[(h_\eps(x))^2\right]=  \E\left[(h^\pm_\eps(x))^2 \right]+\E\left[(h_{\p H}(x))^2 \right]\le \E\left[(h^\pm_\eps(x))^2 \right]+|\log\zeta|,\] for all $x \in H^\pm_\zeta$ and $\eps\in (0,\zeta)$. Note that $|\log\zeta|$ appears from the fact that 
	$$\E\left[ h_{\partial H}^2(x)\right]\le \E\left[h_{\zeta}^2(x)\right]=|\log\zeta|+\log(1-|x|^2)\le |\log\zeta|.$$
	Since $\E\left[(h_\epsilon(x))^2-(h^+_\epsilon(x))^2\right]\leq |\log\zeta|$ for all $\epsilon\in (0,\zeta)$, we have that $\mu_{\gamma h}(\D,\sigma)$ is equal to
\begin{align*}
&	\lim_{\epsilon \to 0}  \left\{ \int_{H\cap \D} e^{\gamma (h_\epsilon^+(x) + h_{\partial D}(x))- \frac{\gamma^2}{2}\E\left[(h_\epsilon(x))^2\right] }dx+
	 \int_{\D\setminus  H} e^{\gamma (h_\epsilon^-(x) + h_{\partial D}(x))- \frac{\gamma^2}{2}\E\left[(h_\epsilon(x))^2\right] }dx\right\}\\
&\hspace{0.1\textwidth}\geq  \zeta^{\gamma^2/2}(\mu_{\gamma h^+}(\D,\sigma \1_{H^+_\zeta}) + \mu_{\gamma h^-}(\D,\sigma \1_{H^-_\zeta}))\inf_{d(z,\partial H)\geq \zeta}e^{\gamma h_{\partial H}(z)}.
	  \end{align*}
	  Here for clarity, we write $\mu_{\gamma h}^{\tilde \sigma}(\D)$ as $\mu_{\gamma h}(\D,\tilde \sigma)$ for $\wt \sigma=\sigma,\sigma \1_{H^+_\zeta}$ and $\sigma \1_{H^-_\zeta}$.
	  
	  Let $K^+$ and $K^-$ be the constants in \eqref{eq:: induction} associated to $\sigma \1_{H^+_\zeta}$ and
	  $\sigma \1_{H^+_\zeta}$ and let $\wt K$ be equal to $\zeta^{-\eta\gamma^22^{n}}K^+K^-$.
Then, by the independence between  $h^+,h^-$ and $h_{\partial H}$,  we have that the expected value of $\exp(-t\mu_{\gamma h}(\D,\sigma))$, conditioned on  $h_{\partial H}$, is upper bounded by
	\begin{align*}
&\E\left[ \exp\left(\ave{-}(t\zeta^{\gamma^2/2}\inf_{d(z,\partial H)\geq \zeta}e^{\gamma h_{\partial H}(z)})(\mu_{\gamma h^+}(\D,\sigma \1_{H^+_\zeta}) + \mu_{\gamma h^-}(\D,\sigma \1_{H^-_\zeta}))\right) \mid h_{\partial H} \right] \\
&\hspace{0.1\textwidth}\le  \widetilde Kt^{-\eta2^{n+1}}  \sup_{d(z,\partial H)\geq \zeta}\exp\left( - \eta \gamma 2^{n+1}h_{\partial H}(z)\right).
	\end{align*}
By Theorem \ref{th.TIS}  and the continuity of $h_{\partial H} $ in  $\overline{ H^\pm_\zeta}$, we have that the expected value of  
$\sup_{d(z,\partial H)\geq \zeta}\exp\left( - \eta \gamma 2^{n+1}h_{\partial H}(z)\right)$ is finite,
which concludes the  proof of \eqref{eq:: induction}.
Now,	we are left with the proof of Claim \ref{claim::separation}.

	\begin{proof}[Proof of Claim \ref{claim::separation}]The fact that $\sigma$ is non-atomic and $\sigma(\D)<\infty$ implies that there are at most countably  many straight lines $\ell$ such that $\sigma(\ell)>0$. Thus, there exists a slope $a\in \R$, such that all straight lines with slope $a$ do not have $\sigma$-mass. WLoG we may assume that $a=0$ satisfies this property. Now, define $E_c$ as follows
		\[E_c:=\D\cap\{z:\Im(z)\geq c\},\] 
		where $\Im(z)$ denotes the imaginary part of $z$.
		Note that $c\mapsto \sigma(E_c)$ is a uniformly continuous function. Thus, there exists a $\bar \zeta>0$ such that for all $c$ the measure of  $\{z:|\Im(z)-c|\leq \bar \zeta\}$ is smaller than or equal to $\sigma(\D)/4$. Furthermore, thanks to the intermediate value theorem, there exists $\bar c$ such that $\sigma(E_{\bar c})=\sigma(\D)/2$. We conclude by taking $H=E_{\bar c}$ and $\zeta=\bar \zeta$. 
	\end{proof}

\bibliographystyle{alpha}

\end{document}